\documentclass[11pt]{article}
\usepackage{cite}
\usepackage{mathrsfs}
\usepackage{amsfonts}
\usepackage{amsmath}
\usepackage{amsfonts,amssymb,color}
\usepackage{dsfont}
\usepackage{curves}
\usepackage{mathrsfs}
\usepackage{pifont}
\usepackage{amssymb}
\usepackage{latexsym,amsmath,amssymb,amsfonts,epsfig,graphicx,cite,psfrag}
\usepackage{eepic,color,colordvi,amscd}

\newtheorem{theorem}{Theorem}[section]

\newtheorem{lemma}{Lemma}[section]
\newtheorem{corollary}{Corollary}[section]

\newtheorem{claim}{Claim}[section]
\newtheorem{conjecture}{Conjecture}[section]

\newcommand{\qed}{\hfill\rule{0.5em}{0.809em}}

\def\emptyset{\mbox{{\rm \O}}}

\def\qed{\hfill \rule{4pt}{7pt}}

\def\pf{\noindent {\it Proof. }}

\begin{document}
	
	\title{Coloring some $(P_6,C_4)$-free graphs with $\Delta-1$ colors}
	\author{Ran Chen$^{1,}$\footnote{Email: 1918549795@qq.com},  \; Di Wu$^{1,}$\footnote{
			Email: 1975335772@qq.com},  \; Xiaowen Zhang$^{2,}$\footnote{Corresponding author. Email: xiaowzhang0128@126.com}\\\\
		\small $^1$Institute of Mathematics, School of Mathematical Sciences\\
		\small Nanjing Normal University, 1 Wenyuan Road,  Nanjing, 210023,  China\\
		\small $^2$School of Mathematical Sciences\\
		\small East China Normal University, Shanghai, 200241, China
	}
	%\small $^2$School of Mathematics, Southeast University, 2 SEU Road, Nanjing, 211189, China}
\date{}

\maketitle
%\;\; and \;\; Yian Xu$^{2,}$\footnote{Email: yian$\_$xu@seu.edu.cn}
%\small  Sophie  Spirkl$^{1,}$\footnote{Email: sophie.spirkl@uwaterloo.ca.},\;\;
%\small $^1$Combinatorics and Optimization, University of Waterloo\\
%\small Waterloo, Ontario, Canada N2L 3G1\\

\begin{abstract}
	The Borodin-Kostochka Conjecture states that for a graph $G$, if $\Delta(G)\geq9$, then $\chi(G)\leq\max\{\Delta(G)-1,\omega(G)\}$. We use $P_t$ and $C_t$ to denote a path and a cycle on $t$ vertices, respectively. Let $C=v_1v_2v_3v_4v_5v_1$ be an induced $C_5$. A {\em $C_5^+$} is a graph obtained from $C$ by adding a $C_3=xyzx$ and a $P_2=t_1t_2$ such that (1) $x$ and $y$ are both exactly adjacent to $v_1,v_2,v_3$ in $V(C)$, $z$ is exactly adjacent to $v_2$ in $V(C)$, $t_1$ is exactly adjacent to $v_4,v_5$ in $V(C)$ and $t_2$ is exactly adjacent to $v_1,v_4,v_5$ in $V(C)$, (2) $t_1$ is exactly adjacent to $z$ in $\{x,y,z\}$ and $t_2$ has no neighbors in $\{x,y,z\}$. In this paper, we show that the Borodin-Kostochka Conjecture holds for ($P_6,C_4,H$)-free graphs, where $H\in \{K_7,C_5^+\}$. This generalizes some results of Gupta and Pradhan in \cite{GP21,GP24}.
	
	\begin{flushleft}
	\end{flushleft}
	
\end{abstract}

\section{Introduction}

All graphs considered in this paper are finite and simple. We use $P_k$ and $C_k$ to denote a {\em path} and a {\em cycle} on $k$ vertices respectively, and follow \cite{BM08} for undefined notations and terminology. Let $G$ be a graph, and let $X$ be a subset of $V(G)$. We use $G[X]$ to denote the subgraph of $G$ induced by $X$, and call $X$ a {\em clique} ({\em stable set}) if $G[X]$ is a complete graph (has no edges). The {\em clique number} $\omega(G)$ of $G$ is the maximum size taken over all cliques of $G$.

Let $G$ and $H$ be two vertex disjoint graphs. The {\em union} $G\cup H$ is the graph with $V(G\cup H)=V(G)\cup V(H)$ and $E(G\cup H)=E(G)\cup E(H)$. The {\em join} $G+H$ is the graph with $V(G+H)=V(G)\cup V(H)$ and $E(G+H)=E(G)\cup E(H)\cup\{xy : x\in V(G), y\in V(H)$$\}$. We say that $G$ induces $H$ if $G$ has an induced subgraph isomorphic to $H$, and say that $G$ is $H$-free otherwise. Analogously, for a family $\cal H$ of graphs, we say that $G$ is ${\cal H}$-free if $G$ induces no member of ${\cal H}$. 

A {\em hole} of $G$ is an induced cycle of length at least 4, and a {\em $k$-hole} is a hole of length $k$. A $k$-hole is called an {\em odd hole} if $k$ is odd, and is called an {\em even hole} otherwise. 

Let $N_G(v)$ be the set of vertices adjacent to $v$ and $d_G(v)=|N_G(v)|$. If it does not cause any confusion, we usually omit the subscript $G$ and simply write $N(v)$ and $d(v)$. Let $\Delta(G)$(resp. $\delta(G)$) denote the maximum(resp. minimum) degree of $G$. 

For positive integer $k$, a $k$-{\em coloring } of $G$ is a function $\phi: V(G)\rightarrow \{1,\cdots,k\}$, such that for each edge $uv$ of $G$, $\phi(u)\ne \phi(v)$. The {\em chromatic number} of $G$, denoted by $\chi(G)$, is the minimum number $k$ for which $G$ has a $k$-coloring. A trivial lower bound of $\chi(G)$ is $\omega(G)$. And by greedy coloring, a trivial upper bound of $\chi(G)$ is $\Delta(G)+1$. In 1941, Brooks \cite{BK77} observed that for a graph $G$, the chromatic number is at most $\Delta(G)$ unless $G$ is a complete graph or an odd hole, in which case these graphs achieve the bound.

\begin{theorem} {\cite{BK77}} 
	For a graph $G$, if $\Delta(G)\geq 3$, then $\chi(G)\leq $ max$\{\Delta(G),\omega(G)\}$.
	
\end{theorem}

In 1977, Borodin and Kostochka conjectured that Brooks' bound can be further improved if
$\Delta(G)\geq 9$.

\begin{conjecture}{\cite{B41}}\label{fredom}
	For a graph $G$, if $\Delta(G)\geq 9$, then $\chi(G)\leq $ max$\{\Delta(G)-1,\omega(G)\}$.
\end{conjecture}

By Brooks' Theorem, each graph $G$ with $\chi(G)>\Delta(G)\geq 9$ contains $K_{\Delta(G)+1}$.
This means Conjecture \ref{fredom} is equivalent to the statement that each $G$ with
$\chi(G)=\Delta(G)\geq 9$ contains $K_{\Delta(G)}$. Note that it was shown in \cite{CLR22} that Conjecture \ref{fredom} cannot be strengthened by making $\Delta(G)\ge 8$ or $\omega(G)\le\Delta(G)-2$.
In 1999, Reed \cite{R99} proved this conjecture for every graph $G$ with $\Delta(G)$ sufficiently large.

\begin{theorem} {\cite{R99}} 
	Every graph with $\chi(G)=\Delta(G)\geq 10^{14}$ contains $K_{\Delta(G)}$.
\end{theorem}

%Before introducing further results about the process of Conjecture \ref{fredom}, we need to define more configurations. A {\em diamond} is a graph obtained from two triangles that share exactly one edge. Let $u_1u_2u_3u_4u_5$ be an induced $P_5$. A {\em kite$^+$} is a graph obtained from $u_1u_2u_3u_4u_5$ by adding a new vertex which is exactly adjacent to $u_1,u_2,u_3$. A {\em flag} is a graph obtained from a $K_4$ by adding a pendant edge, a {\em flag$^+$} is a graph obtained from a flag by subdividing the pendant edge of the flag. A {\em bull} is a graph consisting of a triangle with two disjoint pendant edges. A {\em tripod} is a graph consisting of a triangle with three disjoint pendant edges. A {\em hammer} is a graph obtained by identifying one vertex of a triangle and one end vertex of a $P_3$. A {\em gem} is the graph obtained from a $P_4$ by adding a vertex that is adjacent to $V(P_4)$. A {\em butterfly} is the graph obtained from two trangles that share exactly one vertex. A {\em HVN} is a graph obtained from $K_4$ togethr with one more vertex which is adjacent to exactly two vertices of $K_4$. A {\em house} is a complement of $P_5$. Let $v_1v_2v_3v_4v_5v_1$ be an induced $C_5$. A {\em crown} is the $K_1+K_{1,3}$. 

Let $C=v_1v_2v_3v_4v_5v_1$ be an induced $C_5$. A {\em $C_5^+$} is a graph obtained from $C$ by adding a $C_3=xyzx$ and a $P_2=t_1t_2$ such that (1) $x$ and $y$ are both exactly adjacent to $v_1,v_2,v_3$ in $V(C)$, $z$ is exactly adjacent to $v_2$ in $V(C)$, $t_1$ is exactly adjacent to $v_4,v_5$ in $V(C)$ and $t_2$ is exactly adjacent to $v_1,v_4,v_5$ in $V(C)$, (2) $t_1$ is exactly adjacent to $z$ in $\{x,y,z\}$ and $t_2$ has no neighbors in $\{x,y,z\}$ (See Figure 1 for more forbidden configurations).  Note that $\{v_4,v_5,t_2,v_1,v_2,z\}$ induces a $kite^+$, $\{x,y,v_2,v_1,t_2,v_4\}$ induces a $flag^+$, $\{t_1,v_4,v_5,z,v_3,v_1\}$ induces a $tripod$, $\{x,y,v_1,z,v_3\}$ induces a $crown$, $\{v_1,t_2,v_5,t_1,v_4\}$ induces an $HVN$, $\{v_3,y,v_2,x,v_1\}$ induces a $K_5-e$ and $\{x,v_1,t_2,v_2,v_5\}$ induces a $butterfly$.

\begin{figure}[htbp]\label{fig-1}
	\begin{center}
		\includegraphics[width=11cm]{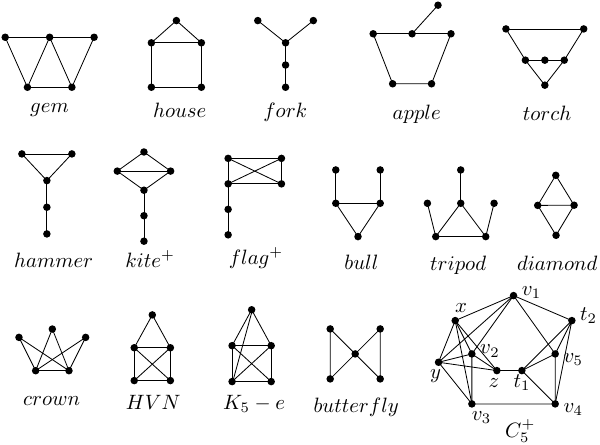}
	\end{center}
	\vskip -15pt
	\caption{Illustration of some forbidden configurations.}
\end{figure}

In 2013, Cranston and Rabern \cite{CR13} proved Conjecture \ref{fredom} for $K_{1,3}$-free graphs. Actually, scholars had been exploring whether Conjecture \ref{fredom} holds on $K_{1,3}$-free graphs for a long time before Cranston and Rabern proved this conclusion. We refer interested readers to \cite{D82,KS86} for this process. Recently, Cranston  {\em et al} \cite{CLR22} proved Conjecture \ref{fredom} for ($P_5, gem$)-free graphs.  Also, Gupta and Pradhan \cite{GP21} proved it for ($P_5, C_4$)-free graphs, Chen {\em et al} \cite{Clz24} proved it for ($P_3\cup P_2, house$)-free graphs. (See Figure 1 for the configurations of $gem$ and $house$.)

In 2021, Dhurandhar \cite{D} proved Conjecture \ref{fredom} for $(P_4\cup K_1)$-free graphs, $P_5$-free graphs and $fork$-free graphs.  This not only generalizes the results of Cranston {\em et al} \cite{CR13,CLR22} and others, as well as Gupta and Pradhan \cite{GP21}, but also provides a new and clever approach for Conjecture \ref{fredom}. Chen {\em et al} proved Conjecture \ref{fredom} holds for $hammer$-free graphs \cite{Cll24} and odd-hole-free graphs \cite{Cllz24}. Wu and Wu \cite{WW23} proved Conjecture \ref{fredom} holds for ($P_6, apple, torch$)-free graphs. (See Figure 1 for the configurations of $fork, hammer, apple$ and $torch$.)

Very recently, Gupta and Pradhan \cite{GP24} proved for each ($P_6,C_4, diamond$)-free graph $G$, if $\Delta(G)\geq5$, then $\chi(G)\leq\max\{\Delta(G)-1,\omega(G)\}$, and they also proved Conjecture \ref{fredom} holds for ($P_6,C_4, H$)-free graphs where $H\in\{bull,K_6\}$ . In this paper, we prove Conjecture \ref{fredom} holds for ($P_6, C_4, H$), where $H\in\{K_7,C_5^+\}$, this is equivalent to the following two theorems. (See Figure 1 for the configurations of $diamond$ and $bull$.)

\begin{theorem}\label{K_7}
	Let $G$ be a $(P_6,C_4,K_7)$-free graph. If $\Delta(G)\geq9$, then $\chi(G)\leq  \max\{\Delta(G)-1,\omega(G)\}$.
\end{theorem}

\begin{theorem} \label{mian theorem}
	Let $G$ be a $(P_6, C_4, C_5^+)$-free graph. If $\Delta(G)\geq 9$, then $\chi(G)\leq  \max\{\Delta(G)-1,\omega(G)\}$.
\end{theorem}

As a corollary, Conjecture \ref{fredom} holds for class of ($P_6,C_4,H$)-free graphs, where $H\in \{kite^+,\\flag^+,tripod,HVN,K_5-e,crown,butterfly\}$. Also, the classes of $(P_5,C_4)$-free graphs,  ($P_6,C_4, \\diamond$)-free graphs and ($P_6,C_4, bull$)-free graphs are all subclasses of $(P_6,C_4,C_5^+)$-free graphs. This generalizes some results of Gupta and Pradhan in \cite{GP21,GP24}.

We will introduce a few more notations, list several useful lemmas, and prove Theorem~\ref{K_7} in section 2. And we will prove Theorem \ref{mian theorem} in Section 3.

\section{Notations and preliminary results}

Let $G$ be a graph, $u,v\in V(G)$, and let $X$ and $Y$ be two subsets of $V(G)$. We simply write $u\sim v$ if $uv\in E(G)$, and write $u\not\sim v$ if $uv\not\in E(G)$.  We say that $v$ is {\em complete} to $X$ if $v$ is adjacent to all vertices of $X$, and say that $v$ is {\em anticomplete} to $X$ if $v$ is not adjacent to any vertex of $X$. We say that $X$ is complete (resp. anticomplete) to $Y$ if each vertex of $X$ is complete (resp. anticomplete) to $Y$.

For a given coloring $\phi$ of $G$, we say $u$ a $m$-vertex if $u\in V(G)$ is colored with $m$. We say $v$ has a $m$-vertex if there exists an $m$-vertex in $N(v)$, in particular, if the $m$-vertex of $v$ is unique, we say $m$ a unique color of $v$ and say the $m$-vertex a unique $m$-vertex of $v$. Moreover, we say $m$ a repeat color of $v$ if there exist more than one $m$-vertex in $N(v)$; and each vertex in $N(v)$ with color $m$ is called a repeat $m$-vertex of $v$. We say $v$ has no missing colors if $N(v)$ received all other colors in $\phi(V(G))$ except the color of $v$. 

We say that a graph $G$ is {\em $k$-vertex-critical} if $G$ has chromatic $k$ and removing any vertex from $G$ results in a graph that has chromatic $(k-1)$. A class of graphs $\mathcal{G}$ is called {\em hereditary} if every induced subgraph of any graph in $\mathcal{G}$ also belongs to $\mathcal{G}$. The following lemmas are very important for our proof.

\begin{lemma}\label{delta}\cite{D51}
	If $G$ is a $k$-vertex-critical graph, then $\delta(G)\geq k-1$. 
\end{lemma}

\begin{lemma}\label{B-K-1}\cite{C76,CLR22,K80}
	Let ${\cal G}$ be a hereditary class of graphs. If Conjecture \ref{fredom} holds for all graphs $G\in {\cal G}$ with $\Delta(G)=9$, then it holds for all graphs in ${\cal G}$.
\end{lemma}

\begin{lemma}\label{B-K-2}\cite{GP21}
	Let ${\cal G}$ be a hereditary class of graphs. If $G\in {\cal G}$ is a counterexample to Conjecture \ref{fredom} with $|V(G)|$ minimal and $\Delta(G)=9$, then $G$ must be vertex-critical.
\end{lemma}

 We say $G$ {\em a relaxed graph} \cite{WW23} if $G$ is a counterexample to Conjecture \ref{fredom} with $|V(G)|$ minimal and $\Delta(G)=9$. Lemma \ref{key} can be directly derived from Lemmas \ref{delta}, \ref{B-K-1} and \ref{B-K-2}.

\begin{lemma}\label{key}\cite{D,Cll24,WW23}
	Let ${\cal G}$ be a hereditary class of graphs, and let $G\in {\cal G}$ be a relaxed graph. Then $G$ is $9$-vertex-critical, $\omega(G)\le8$ and $d(v)=8$ or $9$ for each $v\in V(G)$.
\end{lemma}

Let ${\cal G}$ be a hereditary class of graphs. If Conjecture \ref{fredom} does not hold for ${\cal G}$, then there exists a graph $G\in {\cal G}$ which is a relaxed graph by Lemma~\ref{B-K-1}. By Lemma \ref{key}, we have that $G$ is $9$-vertex-critical and $d(v)=8$ or $d(v)=9$ for each $v\in V(G)$. Choose a $9$-degree vertex $u\in V(G)$ with $N(u)=\{u_1,u_2,\dots,u_7,x,y\}$. Let $\phi(V(G\setminus\{u\}))=\{1,2,\dots,8\}$ be an 8-coloring of $G\setminus\{u\}$ such that $\phi(u_i)=i$ and $\phi(x)=\phi(y)=8$, for $i\in\{1,2,\dots,7\}$. If such $u$ and $\phi$ exist, we say $G$ has a $(u,\phi)$.

For simplicity,  let $[n]=\{1,2,\cdots,n\}$ and $[u_n]=\{u_1,u_2,\cdots,u_n\}$. By Lemma \ref{key}, every relaxed graph $G$ has a $(u,\phi)$ and  $\omega(G)\le8$.

\medskip

Since the proof of Theorem~\ref{K_7} is straightforward by the following lemma, we provide it directly here.

\begin{lemma}\label{P_6,C_4}\cite{KM19}
	Let $G$ be a $(P_6,C_4)$-free graph. Then $\chi(G)\leq\lceil \frac{5}{4}\omega(G)\rceil$.
\end{lemma}

\noindent\textbf{{\em Proof of Theorem~\ref{K_7}} : } Suppose to its contrary that here exists a $(P_6,C_4,K_7)$-free graph $G$ which is a relaxed graph. By Lemma \ref{key}, $G$ is 9-vertex-critical. It contradicts the fact that $\chi(G)\leq8$ by Lemma~\ref{P_6,C_4}. \qed

\medskip

Recall that each relaxed graph $G$ satisfies $\omega(G)\le8$. By Theorem~\ref{K_7}, we have the following corollary.

\begin{corollary}\label{Fin}
	Let ${\cal G}$ be the class of $(P_6,C_4)$-free graphs. If Conjecture \ref{fredom} holds for all graphs $G\in {\cal G}$ with $\omega(G)=7$ and $\omega(G)=8$, then it holds for all graphs in ${\cal G}$.
\end{corollary}

Next, we introduce some more useful lemmas for relaxed graphs. 

\begin{lemma}\label{no missing colors}\cite{D,Cll24}
	Let $G$ be a relaxed graph with a $(u,\phi)$. Then $G$ satisfies that 
	
	\medskip
	
	(1) for $i\in[7]$, $u_i$ has no missing colors $($that is, $\phi(N_{G\setminus\{u\}}(u_i)\cup\{u_i\})=[8]$$)$, and thus $u_i$ has at most one repeat color.
	
	\medskip
	
	(2) for $i,j\in[7]$ and $i\ne j$, if $u_i\not\sim u_j$, then $|N(u_i)\cap N(u_j)\cap[u_7]|\leq2$.
\end{lemma}

For $u,v\in V(G)$, we say an induced path $P$ with endvertices $u$ and $v$ in $G$ an $(i,j)$-$uv$-path if vertices on $P$ are alternately colored with $i$ and $j$. Lemma \ref{path} was essentially proved in \cite{WW23}, we prove it again in a simpler form.

\begin{lemma}\label{path}\cite{WW23}
	Let $G$ be a relaxed graph with a $(u,\phi)$. For $i,j\in[7]$ and $i\ne j$, if $u_i\not\sim u_j$, then there exists an $(i,j)$-$u_iu_j$-path, and $v$ has no missing colors $($that is, $\phi(N_{G\setminus\{u\}}(v)\cup\{v\})=[8]$$)$, where $v$ is any internal vertex of the path.
\end{lemma}

\pf We may by symmetry assume that $u_1\not\sim u_2$ and there does not exist a $(1,2)$-$u_1u_2$-path. Let $H$ denote the subgraph of $G$ induced by the vertices colored by 1 or 2. It is certain that $u_1$ and $u_2$ belong to different components of $H$. Let $H'$ be the component of $H$ which contains $u_1$. So, we can interchange the colors 1 and 2 in $H'$. Now, we can assign $1$ to $u$, and thus $G$ has a proper 8-coloring, a contradiction. 

Moreover, if for each $(1,2)$-$u_1u_2$-path, there exists an internal vertex on the path, say $v$, such that $v$ has a missing color in $\{3,4,\cdots,8\}$, then we can assign the color to $v$. Consequently, there does not exist a $(1,2)$-$u_1u_2$-path, a contradiction.
\qed

\medskip

Recall that $\phi$ is an 8-coloring of $G\setminus \{u\}$ and $\phi(V(G\setminus\{u\}))=[8]$, in the proof of the following lemmas, we do not precolor $u$. 
\begin{lemma}\label{simp}
	Let $G$ be a relaxed graph with a $(u,\phi)$. Suppose $v\in V( G\setminus\{u\})$ such that $v$ has no missing colors. Then the following properties hold.
	
	\medskip
	
	(1) For $r\in[8]$, if $v\in[u_7]$, then $v$ cannot have three $r$-vertices.
	
	\medskip
	
	(2) If $v\not\in N(u)$, then $v$ cannot have three repeat colors and for $r\in[8]$, $v$ cannot have four $r$-vertices.
	
\end{lemma}
\pf Suppose (1) is not true. Then $d(v)\ge 3+6+1=10>\Delta(G)$, a contradiction.
Suppose (2) is not true. Then $d(v)\geq 6+4=10>\Delta(G)$ if $v$ has three repeat colors and $d(v)\geq 4+6=10>\Delta(G)$ if $v$ has four $r$-vertices for some $r\in[8]$, both are contradictions.\qed

\medskip

Lemma \ref{MD} was essentially proved in \cite{D} and the proof sketch also be described in \cite{WW23}. 

\begin{lemma}\label{MD}\cite{D}
	Let $G$ be a relaxed graph with a $(u,\phi)$. Then one of the following does not hold.
	
	\medskip
	
	(\romannumeral 1) $u_i$ is nonadjacent to at most two $u_k$'s for $1\le i,k\le7$.
	
	\medskip
	
	(\romannumeral 2) $|(N(x)\cup N(y))\cap[u_7]|\ge5$.
\end{lemma}
\pf Let $G$ be a relaxed graph with a $(u,\phi)$. Suppose to its contrary that $G$ satisfies (\romannumeral 1) and (\romannumeral 2). Since $G$ is a minimal counterexample for the Borodin-Kostochka conjecture, we have that $\chi(G)=9, \omega(G)\le8$ and $\chi(G\setminus\{u\})=8$. We divide the proof process into two cases depending on $[u_7]$ is a clique or not.

\medskip

\noindent{\bf Case 1} $[u_7]$ is not a clique.

We may assume that $u_1\not\sim u_2$. Then, we will prove that

\begin{equation}\label{eqa-1}
	\mbox{for $i\in\{3,4,\dots,7\}$, $u_i$ cannot be the unique $i$-vertex of both $u_1$ and $u_2$}.
\end{equation}

On the contrary, we may assume that $u_3$ is the unique $3$-vertex of both $u_1$ and $u_2$. By Lemma \ref{no missing colors}(1), $u_3$ cannot have two repeat colors, and thus we suppose $u_1$ is the unique 1-vertex of $u_3$. Now, we can color $u_1,u_2$ by 3, $u_3$ by 1, $u$ by 2. Then we have a proper 8-coloring of $G$, a contradiction. This proves (\ref{eqa-1}).

Next, we prove that 

\begin{equation}\label{eqa-2}
	\mbox{$|N(u_1)\cap N(u_2)\cap N(u)|\le2$}.
\end{equation}
If this case is not true, we may assume that $\{v_1,v_2,v_3\}\subseteq N(u_1)\cap N(u_2)\cap N(u)$. By Lemma \ref{no missing colors}(1), both $u_1$ and $u_2$ have at most one repeat color. Therefore, $u_1$ and $u_2$ have at least one common unique m-vertex in $\{v_1,v_2,v_3\}$, for some $m\in [8]$, say $v_1$. By (\ref{eqa-1}), $v_1\not\in\{v_3,v_4,\dots,v_7\}$, which implies that $v_1\in\{x,y\}$. Without loss of generality, $v_1=x$.
If $x$ has a missing color $r$, then we can color $x$ by $r$, $u_1$ by 8, $u$ by 1, a contradiction. So, $x$ has no missing colors, that is to say, $x$ has at most one repeat color. Without loss of generality, we may assume that $u_1$ is the unique 1-vertex of $x$. Then we can color $u_1,u_2$ by 8, $x$ by 1, $u$ by 2, a contradiction. This proves (\ref{eqa-2}).

By (\romannumeral 1) and $d(u)=9$, we have that $|N(u_1)\cap N(u_2)\cap N(u)|\ge3$, which contradicts to (\ref{eqa-2}).

\medskip

\noindent{\bf Case 2} $[u_7]$ is a clique.

Without loss of generality, we may assume that $|N(x)\cap [u_7]|\geq |N(y)\cap [u_7]|$. Consequently, by (\romannumeral 2), we may assume that $[u_5]\subseteq (N(x)\cup N(y))$  and  $x$ is complete to $[u_3]$. We divide the proof of Case 2 into two subcases: 1) $x$ has no missing colors, and 2) $x$ has a missing color.

{\bf Subase 2.1} $x$ has no missing colors.

Since $[u_7]$ is a clique and $\omega(G)\leq 8$, $x$ has a nonadjacent vertex $u_k$ in $\{u_4,u_5,\cdots,u_7\}$. We have that $x$ has at most one repeat color because $x$ has no missing colors and $d_{G\setminus\{u\}}(x)\leq8$. Moreover, since $u_k$ has at most one repeat color by Lemma~\ref{no missing colors}(1), we have that there exists a vertex in $[u_3]$ which is the unique $m$-vertex of $x$ and $u_k$ for some $m\in [3]$. We may by symmetry assume that $u_1$ is the unique 1-vertex of both $x$ and $u_k$.

Suppose 8 is the unique color of $u_1$. We may color $u_1$ by 8, $x$, $u_k$ by 1, and $u$ by $k$, a contradiction. So, 8 is a repeat color of $u_1$.

Consequently, we prove that
\begin{equation}\label{x}
	\mbox{$y$ is complete to $[u_7]\setminus \{u_k\}$.}
\end{equation}

By Lemma~\ref{no missing colors}(1), $u_1$ has at most one repeat color, and thus $u_k$ is the unique $k$-vertex of $u_1$. We may color $u_1$ by $k$, color $u_k$, $x$ by 1. Now, 8 is a unique color of $u$, and we obtain a new $(u,\phi')$ of $G$ such that $\phi'(u_i)=i$ for $i\in [7]\setminus\{1,k\}$, $\phi'(u_1)=k$, $\phi'(y)=8$ and $\phi'(\{u_k,x\})=1$. For $j\in [7]\setminus\{k\}$, if $y\not\sim u_j$, then we can obtain a contradiction by Case 1. Therefore, $y$ is complete to $[u_7]\setminus \{u_k\}$. This proves (\ref{x}).

By (\ref{x}) and $|N(x)\cap [u_7]|\geq |N(y)\cap [u_7]|$, we have that $x$ is complete to $[u_7]\setminus \{u_k\}$. Obviously, $y$ has at most two repeat colors in $[8]$. Since $\{u_k,x,y\}$ is complete to $[u_7]\setminus \{u_k\}$ and both $u_k$ and $x$ have at most one repeat color, we have that there exists a vertex $u_t\in[u_7]\setminus \{u_k\}$ such that $u_t$ is the unique $t$-vertex of $u_k$, $x$ and $y$. Since $u_t$ has at most one repeat color by Lemma~\ref{no missing colors}(1), we have that $u_k$ is the unique $k$-vertex of $u_t$. Now, we can color $u_t$ by $k$, $u_k$, $x$, $y$ by $t$, and $u$ by 8. Now, $G$ has a proper 8-coloring, a contradiction.

\medskip

{\bf Subcase 2.2} $x$ has a missing color.

We may assume that color $k$ is a missing color of $x$, then $x\not\sim u_k$. 

Suppose $y$ has a missing color, say $l$. Then, we can color $y$ by $l$, $x$ by $k$, and $u$ with 8, a contradiction. So, $y$ has no missing colors, and thus $y$ has at most one repeat color.

Next, we prove that

\begin{equation}\label{y}
	\mbox{$y$ is complete to $[u_7]\setminus \{u_k\}$.}
\end{equation}

We may color $x$ by $k$. Now, 8 is a unique color of $u$, and we obtain a new $(u,\phi'')$ of $G$ such that $\phi''(u_i)=i$ for $i\in [7]\setminus\{k\}$, $\phi''(y)=8$ and $\phi''(\{u_k,x\})=k$. For $j\in [7]\setminus\{k\}$, if $y\not\sim u_j$, then we can obtain a contradiction by Case 1. This proves (\ref{y}).

Since $\{u_k,y\}$ is complete to $[u_7]\setminus \{u_k\}$ and both $u_k$ and $y$ have at most one repeat color by (\ref{y}) and Lemma~\ref{no missing colors}(1), we have that there exists a vertex $u_t\in[u_7]\setminus \{u_k\}$ such that $u_t$ is the unique $t$-vertex of $u_k$ and $y$. 

Since $u_t$ has at most one repeat by Lemma~\ref{no missing colors}(1), we have that either 8 or $k$ is a unique color of $u_t$. If 8 is a unique color of $u_t$, then we can color $u_t$ by 8, $y$, $u_k$ by $t$, and $u$ by $k$, then we obtain a proper 8-coloring of $G$. If $k$ is a unique color of $u_t$, we can color $u_t$, $x$ by $k$, $u_k$, $y$ by $t$, and $u$ by 8, then we obtain a proper 8-coloring of $G$. Both are contradictions. \qed

\medskip

\section{Proof of Theorem~\ref{mian theorem}}

In this section, we will prove Theorem \ref{mian theorem} that Conjecture \ref{fredom} holds for ($P_6, C_4$, $C_5^+$)-free graphs. The following lemma is useful for our proof.

\begin{lemma}\label{main1}
	Let $G$ be a relaxed graph with a $(u,\phi)$. If $G$ is $(P_6, C_4)$-free, then there exists $i\in[7]$ such that $|N(u_i)\cap [u_7]|\leq 3$.
\end{lemma}

\pf
Suppose to its contrary that $|N(u_i)\cap [u_7]|\geq
4$ for $i\in[7]$.

We first prove that
\begin{equation}\label{clique}
	\mbox{$[u_7]\cup\{u\}$ is a clique.}
\end{equation}

Suppose there exist $u_i\not\sim u_j$, where $i,j\in[7]$. Then, we have that $|N(u_i)\cap N(u_j)\cap [u_7]|\geq3$, which contradicts to Lemma \ref{no missing colors}(2). Therefore, $[u_7]\cup\{u\}$ is a clique. This proves (\ref{clique}).

Next, we prove that 
\begin{claim}\label{5}
	$|(N(x)\cup N(y))\cap[u_7]|\ge5.$
\end{claim}
\pf
Suppose to its contrary that $|(N(x)\cup N(y))\cap[u_7]|\leq4$. Without loss of generality, we may assume that $\{x,y\}$ is anticomplete to $\{u_1,u_2,u_3\}$.

First, we prove that
\begin{equation}\label{colors}
	\mbox{both $x$ and $y$ have no missing colors. }
\end{equation}

Suppose each vertex of  $\{x,y\}$ has a missing color. Then, we can assign $x$ and $y$ with the missing color of each of them, and assign 8 to $u$. Consequently, $G$ has a proper 8-coloring, a contradiction. 

Suppose exactly one of $\{x,y\}$ has no missing colors. We may by symmetry assume that $y$ has a missing color $k\in[7]$. Then, we can color $y$ with $k$ to form a new coloring $\phi'$ of $G\setminus\{u\}$ such that $G$ has a $(u,\phi')$ and $k$ is a repeat color of $u$.  Now, we have that for $t\in ([u_7]\setminus\{u_k\})\cup\{x\}$, $|N(t)\cap (([u_7]\setminus\{u_k\})\cup\{x\})|\geq4$. Let $t_1,t_2\in ([u_7]\setminus\{u_k\})\cup\{x\}$. If $t_1\not\sim t_2$, then $|N(t_1)\cap N(t_2)\cap ([u_7]\setminus\{u_k\})\cup\{x\}|\geq3$. But by Lemma~\ref{no missing colors}(2), $|N(t_1)\cap N(t_2)\cap ([u_7]\setminus\{u_k\})\cup\{x\}|\leq2$, a contradiction. Therefore, $t_1\sim t_2$ and thus ($[u_7]\setminus\{u_k\})\cup\{x\}$ is a clique. And now, We have that $N(x)\cap [u_3]\ne\emptyset$, a contradiction. This proves (\ref{colors}).

Next, we prove that 
\begin{equation}\label{path2}
	\mbox{for $i\in[7]$, there exists an $(i,8)$-$u_ix$-path or $(i,8)$-$u_iy$-path.}
\end{equation}

Suppose to its contrary that there does not exist an $(i,8)$-$u_ix$-path or $(i,8)$-$u_iy$-path for some $i\in[7]$. Then, let $H$ denote the subgraph of $G$ which induced by the vertices which is colored by $i$ or 8. Let $H'$ be the component of $H$ which contains $u_i$. It is certain that $x\not\in V(H')$ and $y\not\in V(H')$.  Now, we may interchange the colors 8 and $i$ in $H'$. We can assign $i$ to $u$, and thus $G$ has a proper 8-coloring, a contradiction. This proves (\ref{path2}).

Let $P$ be an $(i,8)$-$u_ix$-path or $(i,8)$-$u_iy$-path and $v$ is an internal vertex of $P$. If $v$ has a missing color, then we can assign the color to $v$. So, by (\ref{path2}), we may assume that each internal vertex of $P$ has no missing colors. 

By (\ref{path2}), we have that the length of each $(i,8)$-$u_ix$-path (or $(i,8)$-$u_iy$-path) is odd. Moreover, since $\{x,y\}$ is anticomplete to $\{u_1, u_2, u_3\}$ and $G$ is $P_6$-free, without loss of generality, let $P_a=u_1v_8v_1a$, $P_b=u_2v_8'v_2b$ and $P_c=u_3v_8''v_3c$ be a $(1,8)$-$u_1a$-path, a $(2,8)$-$u_1b$-path and a $(3,8)$-$u_1c$-path respectively, where $a,b,c\in\{x,y\}$, such that all of $\{v_8,v_8',v_8'',v_1,v_2,v_3\}$ have no missing colors (it is possible that $|\{v_8,v_8',v_8''\}|\ne3$). Without loss of generality, we set $a=x$.

Now, without loss of generality, we may divide the proof of Claim~\ref{5} into three cases: 1) $v_8=v_8'=v_8''$, 2) $|\{v_8,v_8',v_8''\}|=3$, and 3) $|\{v_8,v_8',v_8''\}|=2$.

\medskip

{\bf Case 1} $v_8=v_8'=v_8''$.

Now, $v_8$ has two 1-vertices, two 2-vertices and two 3-vertices. It contradicts Lemma \ref{simp}(2), a contradiction. 

\medskip

{\bf Case 2} $|\{v_8,v_8',v_8''\}|=3$.

By (\ref{clique}), $[u_7]\cup \{u\}$ is a clique. Note that $x\not\sim y$. If $v_1\sim y$, then $\{v_1,y,u,x\}$ induces a $C_4$. This contradiction implies that $v_1\not\sim y$. Since $y$ has no missing colors, we may assume that $w_1$ is a 1-vertex of $y$. It is certain that $w_1\notin\{v_1,u_1\}$. To avoid an induced $C_4$ on $\{w_1,u,x,y\}$, we have that $w_1\not\sim x$. Consequently, we have that $w_1\sim v_8$ as otherwise there exists an induced $P_6=w_1yuu_1v_8v_1$. Now, the path $u_1v_8w_1y$ is a $(1,8)$-$u_1y$-path. Recall that the path $u_1v_8v_1x$ is a (1,8)-$u_1x$ path.

Next, we prove that 

\begin{equation}\label{v_8}
	\mbox{$v_8\not\sim u_2$.}
\end{equation}

Suppose to its contrary that $v_8\sim u_2$. Since $b\in\{x,y\}$ and $G[\{v_2,u,x,y\}]$ is not an induced $C_4$, we may by symmetry assume that $v_2\sim x$ and $v_2\not\sim y$.

If $u_2\sim w_1$, then $\{y,w_1,u_2,u\}$ induces a $C_4$, this contradiction implies that $u_2\not\sim w_1$. To forbid an induced $C_4$ on $\{w_1,v_8,u_2,v_8'\}$, we have that $w_1\not\sim v_8'$. Moreover, To forbid an induced $C_4$ on $\{v_2,v_8,v_8',u_2\}$, we have that $v_2\not\sim v_8$. So, we can deduce that $w_1\sim v_2$ as otherwise there exists an induced $P_6=yw_1v_8u_2v_8'v_2$. 

To forbid an induced $C_4$ on $\{w_1,v_1,v_2,v_8\}$ or $\{u,u_2,x,v_1\}$ or $\{u,u_3,x,v_1\}$ or $\{u,u_3,x,v_2\}$, we have that $v_1$ is anticomplete to $\{v_2,u_2,u_3\}$ and $u_3\not\sim v_2$. Therefore, we have that $u_3\sim v_8$ as otherwise there exists an induced $P_6=u_3u_2v_8v_1xv_2$. 

Now, $u_3$ is complete to $\{v_8,v_8''\}$. By Lemma \ref{simp}(1), we have that $u_3\not\sim v_8'$. Consequently, to avoid an induced $C_4=uyw_1u_3u$, we have that $w_1\not\sim u_3$. But then, there exists an induced $P_6=yw_1v_2v_8'u_2u_3$, a contradiction. This proves (\ref{v_8}).

Consequently, we prove that

\begin{equation}\label{v_8'}
	v_8'\not\sim u_1.
\end{equation}

Since $b\in\{x,y\}$ and $G[\{v_2,u,x,y\}]$ is not an induced $C_4$, we have that $v_2$ is adjacent exactly to one vertex in $\{x,y\}$. 

If $v_2\not\sim x$, then let $w_2\notin\{u_2,v_2\}$ be a 2-vertex of $x$. We have that $w_2\sim v_8'$ as otherwise there is an induced $P_6=w_2xuu_2v_8'v_2$. Now, $u_2v_8'w_2x$ is a $(2,8)$-$u_2x$-path and $u_2v_8'v_2y$ is a $(2,8)$-$u_2y$-path. 

If $v_2\not\sim y$, then let $t_2\notin\{u_2,v_2\}$ be a 2-vertex of $y$. We have that $t_2\sim v_8'$ as otherwise there is an induced $P_6=t_2yuu_2v_8'v_2$. Now, $u_2v_8'v_2x$ is a $(2,8)$-$u_2x$-path and $u_2v_8't_2y$ is a $(2,8)$-$u_2y$-path. 

Therefore, with the same arguments as in the proof of (\ref{v_8}), we have that $v_8'\not\sim u_1$. This proves (\ref{v_8'}).

To avoid an induced $C_4=uyw_1u_2u$, it follows that $w_1\not\sim u_2$. Then $w_1\sim v_8'$ as otherwise there exists an induced $P_6=yw_1v_8u_1u_2v_8'$, a contradiction. Furthermore, if $v_1\not\sim v_8'$, then there exists an induced $P_6=v_1xuyw_1v_8'$. So,  $v_1\sim v_8'$. But now, there exists an induced $C_4=w_1v_8v_1v_8'w_1$, a contradiction. 

\medskip

{\bf Case 3} $|\{v_8,v_8',v_8''\}|=2$.

We may by symmetry assume that $v_8\ne v_8'=v_8''$. To avoid an induced $C_4$ on $\{u,x,y,v_2\}$, we have that $v_2$ is not complete to $\{x,y\}$. Without loss of generality, we may assume that $v_2\not\sim y$. 

Since $y$ has no missing colors by (\ref{colors}), let $w_2$ be a 2-vertex of $y$. Also, since $v_8'$ is complete to $\{v_2,v_3,u_2,u_3\}$ and $v_8'$ has no missing colors, we have that $v_8'\not\sim w_2$. But now, there exists an induced $P_6=w_2yuu_2v_8'v_2$, a contradiction. 

This completes the proof of Claim~\ref{5}.\qed

\medskip

By Claim~\ref{5}, we have that $|(N(x)\cup N(y))\cap[u_7]|\ge5$. By Lemma~\ref{MD}, there exists $i\in[7]$ such that $|N(u_i)\cap [u_7]|\leq 3$, a contradiction. It completes the proof of Lemma~\ref{main1}. \qed

\medskip

\noindent\textbf{{\em Proof of Theorem~\ref{mian theorem}} : } Suppose that Theorem~\ref{mian theorem} does not hold. Let $G$ be a relaxed graph with a $(u,\phi)$, and $G$ is $(P_6, C_4, C_5^+)$-free. By Lemma~\ref{main1}, without loss of generality, we may assume that $u_1$ is anticomplete to $\{u_2,u_3,u_4\}$. By Lemma~\ref{path}, there must exist a $(1,2)$-$u_1u_2$-path in $G$, say $P^{1,2}$ such that any internal vertex on the path has no missing colors. Since $G$ is $P_6$-free, we may assume that $P^{1,2}=u_1v_2v_1u_2$. Similarly, let $P^{1,3}=u_1v_3v_1'u_3$ be a $(1,3)$-$u_1u_3$-path in $G$ and $P^{1,4}=u_1v_4v_1''u_4$  be a $(1,4)$-$u_1u_4$-path in $G$ such that all of $\{v_3,v_1',v_4,v_1''\}$ have no missing colors (it is possible that $|\{v_1,v_1',v_1''\}|\ne3$). 

Now, we divide the proof of Theorem~\ref{mian theorem} into three cases: 1) $v_1=v_1'=v_1''$, 2) $|\{v_1,v_1',v_1''\}|=3$, and 3) $|\{v_1,v_1',v_1''\}|=2$.

\medskip

{\bf Case 1} $v_1=v_1'=v_1''$.

Now, $v_1$ has two 2-vertices, two 3-vertices and two 4-vertices, this contradicts Lemma \ref{simp}(2).

\medskip

{\bf Case 2} $|\{v_1,v_1',v_1''\}|=3$.

By Lemma~\ref{no missing colors}(1), $u_1$ has at most one repeat color in $\{2,3,4\}$. We may by symmetry suppose that 2 is a unique color of $u_1$. Notice that $v_1'$ has no missing colors, and so we next consider the 2-vertex of $v_1'$.

We first prove that

\begin{equation}\label{u_2}
	\mbox{$u_2\not\sim v_1'$.}
\end{equation} 

Suppose to its contrary that $u_2\sim v_1'$. To avoid an induced $C_4$ on $\{u,u_1,u_2,v_4\}$, we have that $v_4\not\sim u_2$. By Lemma~\ref{no missing colors}(1), both $u_2$ and $v_4$ have no missing colors. If  $u_2\sim v_1''$, then $u_2$ has three 1-vertices, which contradicts Lemma \ref{simp}(1). If $v_4$ is complete to $\{v_1,v_1'\}$, then $v_4$ has four 1-vertices, which contradicts Lemma \ref{simp}(2). So, $u_2\not\sim v_1''$ and $v_4$ is not complete to $\{v_1,v_1'\}$. But then, $\{v_1,u_2,u,u_1,v_4,v_1''\}$ or $\{v_1',u_2,u,u_1,v_4,v_1''\}$ induces a $P_6$, a contradiction. This proves (\ref{u_2}).

Next, we prove that
\begin{equation}\label{v_2}
	\mbox{$v_2\not\sim v_1'$.}
\end{equation} 

Suppose to its contrary that $v_2\sim v_1'$. By Lemma \ref{simp}(2), we have that $v_2\not\sim v_1''$. Suppose $u_2\sim v_1''$. By (\ref{u_2}), $u_2\not\sim v_1'$. But now, there exists an induced $P_6=v_1'v_2u_1uu_2v_1''$, a contradiction. So, $u_2\not\sim v_1''$.

To avoid an induced $C_4=uu_1v_4u_2u$, we have that $v_4\not\sim u_2$. We can deduce that $v_1\sim v_4$ as otherwise there exists an induced $P_6=v_1u_2uu_1v_4v_1''$. Now, $v_4$ has three 1-vertices, and so $v_4\not\sim v_1'$ by Lemma \ref{simp}(2). To avoid an induced $C_4$ on $\{u,u_1,u_3,v_4\}$, $v_4\not\sim u_3$. If $v_1\not\sim u_3$, then there exists an induced $P_6=v_1'u_3uu_1v_4v_1$, a contradiction. So, $v_1\sim u_3$. To avoid an induced $C_4=uu_1v_2u_3u$, $v_2\not\sim u_3$. But now, $\{v_2,v_1,u_3,v_1'\}$ induces a $C_4$, a contradiction. This proves (\ref{v_2}).

By (\ref{u_2}) and (\ref{v_2}), we have that both $u_2$ and $v_2$ are not the 2-vertex of $v_1'$. Therefore, let $w_2$ be a 2-vertex of $v_1'$. Then $w_2\not\in\{u_2,v_2\}$. Since $v_2$ is a unique 2-vertex of $u_1$ by our assumption, it follows that $u_1\not\sim w_2$. To avoid an induced $C_4=uu_1v_2u_3u$, we have that $v_2\not\sim u_3$. Then $v_1\sim u_3$ as otherwise there exists an induced $P_6=v_1v_2u_1uu_3v_1'$. To avoid an induced $C_4=v_1u_2uu_3v_1$, we have that $u_2\sim u_3$. Now, $u_3$ has at least two 1-vertices $v_1$ and $v_1'$. Since $u_2$ is a 2-vertex of $u_3$, we have that $w_2\not\sim u_3$ by Lemma~\ref{no missing colors}(1). Then, there exists an induced $P_6=w_2v_1'u_3uu_1v_2$, a contradiction. 

\medskip

{\bf Case 3} $|\{v_1,v_1',v_1''\}|=2$.

Without loss of generality, we may assume that $v_1\not=v_1'=v_1''$. By Lemma~\ref{no missing colors}(1), there exists at most one repeat color of $u_1$ in $\{3,4\}$. Without loss of generality, we may assume that 3 is a unique color of $u_1$.

To avoid an induced $C_4$ on $\{u,u_3,v_1',u_4\}$ or $\{u_1,v_3,v_1',v_4\}$, we have that $u_3\sim u_4$ and $v_3\sim v_4$. Also, to forbid a $C_4$ on $\{u,u_1,v_3,u_2\}$, $\{u,u_1v_4,u_2\}$, $\{u,u_1,u_4,v_3\}$ or $\{u,u_1,u_3,v_4\}$, we have that $v_3$ is anticomplete to $\{u_2,u_4\}$ and $v_4$ is anticomplete to $\{u_2,u_3\}$. And, we can deduce that $v_2$ is anticomplete to $\{u_3,u_4\}$,  otherwise there exists an induced $C_4$ on $\{u,u_1,v_2,u_3\}$ or $\{u,u_1,v_2,u_4\}$. We first prove that

\begin{equation}\label{v_1'}
	\mbox{$v_1'$ is adjacent to exactly one vertex in $\{v_2,u_2\}$}.
\end{equation}

Suppose to its contrary. If $v_1'$ is complete to $\{v_2,u_2\}$, then $\{v_2,v_1,u_2,v_1'\}$ induces a $C_4$. So,  $v_1'$ is anticomplete to $\{v_2,u_2\}$. If $u_2\not\sim u_3$, then $v_1\not\sim u_3$ to avoid an induced $C_4=uu_2v_1u_3u$. But now, an induced $P_6=v_2v_1u_2uu_3v_1'$ appears, a contradiction. So, $u_2\sim u_3$. To forbid an induced $P_6=u_1v_2v_1u_2u_3v_1'$, we have that $v_1\sim u_3$.

If $v_2\not\sim v_3$, then $v_1\not\sim v_3$ to avoid an induced $C_4=u_1v_2v_1v_3u_1$, but now there exists an induced $P_6=u_2v_1v_2u_1v_3v_1'$. So, $v_2\sim v_3$. Then $v_1\sim v_3$ as otherwise there exists an induced $P_6=uu_2v_1v_2v_3v_1'$. But now, an induced $C_4=v_1v_3v_1'u_3v_1$ appears, a contradiction. This proves (\ref{v_1'}).

Next, we prove that

\begin{equation}\label{v_1'+}
	\mbox{$v_1'\sim v_2$ and $v_1'\not\sim u_2$.}
\end{equation}

Suppose to its contrary that $v_1'\sim u_2$ and $v_1'\not\sim v_2$ by (\ref{v_1'}). We can deduce that $u_2$ is complete to $\{u_3, u_4\}$ because both $uu_2v_1'u_3u$ and $uu_2v_1'u_4u$ are not induced $C_4$. Since both $v_1v_2u_1uu_3v_1'$ and $v_1v_2u_1uu_4v_1'$ are not induced $P_6$, $v_1$ is complete to $\{u_3,u_4\}$. Also, to avoid an induced $C_4$ on $\{v_1,u_2,v_1',v_3\}$ or $\{v_1,u_2,v_1',v_4\}$, we have that $v_1$ is anticomplete to $\{v_3,v_4\}$. 

Suppose $v_2\not\sim v_3$. Let $w_3\notin\{v_3,u_3\}$ be a 3-vertex of $v_2$ because $v_2$ has no missing colors. Since 3 is a unique color of $u_1$ and $u_1\sim v_3$, we have that $u_1\not\sim w_3$. Since $v_1'$ has no missing colors and $v_1'$ is complete to $\{v_3,v_4,u_3,u_4\}$, we have that $v_1'\not\sim w_3$ as otherwise $d(v_1')>5+5=10$. But now, there is an induced $P_6=w_3v_2u_1v_3v_1'u_3$, a contradiction. Therefore, $v_2\sim v_3$.

Suppose $v_2\not\sim v_4$. Let $w_4\notin\{v_4,u_4\}$ be a 4-vertex of $v_2$ and $w_2\notin\{v_2,u_2\}$ be a 2-vertex of $v_4$ as both $v_2$ and $v_4$ have no missing colors.  Since $v_1'$ has no missing colors and $v_1'$ is complete to $\{v_3,v_4,u_3,u_4\}$, we have that $v_1'\not\sim w_4$ as otherwise $d(v_1')>5+5=10$. By Lemma~\ref{simp}(2), $v_1'\not\sim w_2$. Since $w_4v_2u_1v_4v_1'u_4$ is not an induced $P_6$, we have that $u_1\sim w_4$. By Lemma~\ref{no missing colors}(1) and $u_1$ is complete to $\{v_2,u_4,w_4\}$, we have that 2 is a unique vertex of $u_1$, and so $u_1\not\sim w_2$. To avoid an induced $P_6=w_2v_4u_1uu_2v_1$, $w_2\sim v_1$. Now, we have that $v_1$ is complete to $\{v_2,u_2,w_2,u_4\}$. Since $v_1$ has no missing colors, we have that $v_1\not\sim w_4$ as otherwise $d(v_1)>5+5=10$. To avoid an induced $C_4=w_4v_2v_1w_2w_4$, $w_2\not\sim w_4$. But now, there exists an induced $P_6=w_4v_2v_1w_2v_4v_1'$, a contradiction. Therefore, $v_2\sim v_4$.

Now, $uu_1v_3v_1'u_3u$ is an induced $C_5$, and $G[\{u,u_1,v_3,v_1',u_3,u_2,u_4,v_1,v_2,v_4\}]$ is an induced $C_5^+$, a contradiction. This proves (\ref{v_1'+}).

By (\ref{v_1'+}), $v_1'\sim v_2$ and $v_1'\not\sim u_2$. To avoid an induced $C_4=u_1v_2v_1'v_3u_1$ or $u_1v_2v_1'v_4u_1$, we have that $v_2$ is complete to $\{v_3,v_4\}$. Since both $v_1u_2uu_1v_3v_1'$ and $v_1u_2uu_1v_4v_1'$ are not induced $P_6$, $v_1$ is complete to $\{v_3,v_4\}$. We can also deduce that $v_1$ is anticomplete to $\{u_3,u_4\}$ as otherwise there exists an induced $C_4=v_3v_1u_3v_1'v_3$ or $v_3v_1u_4v_1'v_3$. 

Suppose $u_2\not\sim u_3$. Let $t_3\notin\{v_3,u_3\}$ be a 3-vertex of $u_2$ as $u_2$ has no missing colors by Lemma~\ref{no missing colors}(1). Since $v_1'$ has no missing colors and $v_1'$ is complete to $\{v_3,v_4,u_3,u_4\}$, we have that $v_1'\not\sim t_3$ as otherwise $d(v_1')>5+5=10$. But now, $t_3u_2uu_3v_1'v_3$ is an induced $P_6$, a contradiction. So, $u_2\sim u_3$. Similarly, $u_2\sim u_4$. 

Now, $uu_1v_3v_1'u_3u$ is an induced $C_5$, and $G[\{u,u_1,v_3,v_1',u_3,v_2,v_4,v_1,u_2,u_4\}]$ is an induced $C_5^+$, a contradiction. Therefore $G$ does not exist and completes the proof of Theorem~\ref{mian theorem}.\qed

\section*{Declarations}

\begin{itemize}
	\item \textbf{Funding}\quad This work was  supported by NSFC 11931006.
	\item \textbf{Conflict of interest}\quad The authors declare that they have no known competing financial interests or personal
	relationships that could have appeared to influence the work reported in this paper.
	\item \textbf{Data availibility statement}\quad This manuscript has no associated data.
\end{itemize}

\end{document}